\documentclass[10pt,leqno]{amsart}
\usepackage{amssymb,amsmath,amsthm,amssymb,dsfont}
\usepackage[margin=0.55in]{geometry}
\usepackage{amsmath,amssymb,booktabs,longtable,array,xcolor}
\usepackage{color}
\usepackage[normalem]{ulem}
\usepackage[bookmarksopen,bookmarksdepth=3,colorlinks,citecolor=red,pagebackref,hypertexnames=true]{hyperref}
\usepackage{amsmath,amssymb,booktabs,longtable,array,xcolor}
\usepackage[msc-links,nobysame,non-sorted-cites, initials]{amsrefs}
\usepackage[inline]{enumitem}
\usepackage{mathtools}
\mathtoolsset{showonlyrefs}

\definecolor{darkblue}{RGB}{0,0,160}

\usepackage[lining]{libertine}   
\usepackage{cabin}
\usepackage[libertine]{newtxmath}
\usepackage[utf8]{inputenc}
\usepackage[T1]{fontenc}
\def\eps{\varepsilon}
\def\d{{\rm d}}
\def\R {\mathbb{R}}
\def\N {\mathbb{N}}

\def\Z {{\mathbb Z}}

\newcounter{counter}

\numberwithin{equation}{section}
\numberwithin{counter2}{section}
\newtheorem{proposition}[subsection]{Proposition}
\newtheorem{theorem}[counter]{Theorem}
\newtheorem{corollary}[counter]{Corollary}
\newtheorem{lemma}[subsection]{Lemma}

\theoremstyle{definition}
\newtheorem{definition}[subsection]{Definition}
\newtheorem*{remark*}{Remark}
\newtheorem*{warn*}{A word of warning}

\theoremstyle{plain}
\newcommand{\TT}{\mathbb{ T} }

\newcommand{\Cramer}{\operatorname{Cram\acute{e}r}}

\newcommand{\ind}[1]{\mathds{1}_{{#1}}}
\newcommand{\floor}[1]{\lfloor #1 \rfloor}

\newcommand{\ceiling}[1]{\left\lceil #1 \right\rceil}
\newtheorem{problem}[equation]{Problem}
\usepackage{scalerel} 
\newcommand{\EE}{\mathbb{E}}

\title{Pointwise Ergodic Theorems Along $\lfloor p^c\rfloor$ for Every $c>1$
}

\author[A. Fragkos]{Anastasios Fragkos}
\begin{document}
	\begin{abstract}
	    Let $c>1$ be noninteger, and let $\Lambda$ denote the von Mangoldt function. The averages
\[
A_N^{c} f(x):=\frac1N\sum_{n=1}^{N}
\Lambda(n)f\bigl(T^{\lfloor n^c\rfloor}x\bigr)
\]
satisfy the variational estimate
\[
\bigl\|\mathcal{V}^r(A_N^{c} f:N\in\mathbb{N})\bigr\|_{L^p(X)}
\lesssim_{c,p,r}\|f\|_{L^p(X)}
\]
for every $1<p<\infty$ and $r>2$, on every
$\sigma$-finite measure-preserving system. In particular, the
normalized averages along fractional powers of primes converge
almost everywhere. The main arithmetic ingredient is a uniform bound on the correlations of $\Lambda-1$ with fractional-power exponential phases, yielding arbitrarily  large powers of logarithmic gain. We obtain this estimate through
short-interval polynomial approximation of the fractional phase and higher-order discorrelation estimates. The same approach yields
almost-everywhere convergence to for the
M\"obius-weighted averages along $\left( \floor{n^c} \right)_{n \in \N}$ throughout $1<p<\infty$.
	\end{abstract}
	\maketitle
	\section{Introduction}
	
	Birkhoff's theorem \cite{Birkhoff1931} gives almost everywhere
	convergence of the ordinary ergodic averages for integrable
	functions. The corresponding problem for sparse arithmetic times
	was developed in Bourgain's work
	\cite{Bourgain1988Maximal,Bourgain1988Lp,Bourgain1989}. He established
	pointwise convergence along the squares and, more generally, along
	integer polynomial sequences, first for $L^2$ functions and then
	on $L^p$ for all $p>1$. Thus, on an invertible
	measure-preserving system, the averages
	\[\frac{1}{N}\sum_{n\le N}f(T^{P(n)}x)\] converge almost everywhere whenever
	$P\in\mathbb Z[t]$ and $f\in L^p(X)$, for $p$ in the reflexive range.

	For averages sampled along prime orbits, Wierdl \cite{Wierdl1988} extended
	Bourgain's initial results to $L^p$ for $p \in (1,\infty)$, including
	fixed positive integer powers of primes. Nair \cite{Nair1993}
	subsequently proved pointwise convergence along polynomial images
	of the primes. Quantitative estimates in the prime setting were
	obtained by Mirek, Trojan and Zorin-Kranich \cite{MTZK2017}, who
	established $L^p$ bounds for the $r$-variation,
	for $ 1<p<\infty$ and $r>2$.
	
	We consider averages along fractional powers of primes. Let
	$(X,\nu,T)$ be a $\sigma$-finite
	measure-preserving system and let $c>1$ be noninteger. Our point of focus will be the averages
	\[
	A_N^{c} f(x):=\frac1N\sum_{n=1}^N
	\Lambda(n)f\bigl(T^{\lfloor n^c\rfloor}x\bigr),
	\]
	where $\Lambda$ denotes the von Mangoldt function. The exponent
	$c$ is fixed throughout.
	
	The unweighted counterpart of $A_N^{c}$ is denoted by
	\[
	B_N^{c} f(x):=\frac1N\sum_{n=1}^N f(T^{\lfloor n^c\rfloor}x).
	\]

	The unweighted averages have been studied since Wierdl's
	doctoral thesis \cite{Wierdl1989Thesis}. In \cite[Theorem~4.3]{Wierdl1989Thesis}, he proved
	their almost everywhere convergence for every noninteger $c>1$
	and every $f\in L^p(X)$, for $p>1$. Boshernitzan and Wierdl
	\cite{BW1996}, and Boshernitzan, Kolesnik, Quas and Wierdl
	\cite{BKQW2005}, subsequently studied averages along integer parts
	of Hardy field functions. Their results include many combinations
	of powers and logarithms and relate the averaging properties to
	the growth of the function and its distance from polynomials.
	
	At the $L^1$ endpoint, Urban and Zienkiewicz \cite{UZ2007}
	established weak $(1,1)$ maximal estimates and pointwise
	convergence for the unweighted averages when $1<c<1001/1000$.
	Mirek \cite{Mirek2015Rough} extended the range to $1<c<30/29$
	and treated suitable regularly varying perturbations of $n^c$.
	More recently, Krause and Sun \cite{KrauseSun2025} obtained
	quantitative convergence results for $1<c<7/6$. These endpoint
	results concern the averages $B_N^{c}$.

    The times $\left( \floor{p^c} \right)_{p \in \mathbb{P}}$ exhibit two sources of sparsity: the restriction to primes and their image under $x \mapsto x^c, \; \; c>1.$  Fractional powers of primes have been studied in the scope of recurrence, equidistribution, joint ergodicity and $L^2$ convergence of multiple ergodic averages, see for example \cite{Bergelson2014,BKS2019,Frantzikinakis2022Joint,KoutsogiannisTsinas2026}.   
	
    In addition to that, Bahnson, Daskalakis, Dohadwala
and Shah \cite{BDDS24} proved pointwise convergence throughout
the range $1<c<4/3$. More generally, on a
$\sigma$-finite measure-preserving system, they considered
\[
P_N^h f(x):=\frac1{\pi(N)}
\sum_{\substack{q\le N\\q\in\mathbb{P}}}
f\bigl(T^{\lfloor h(q)\rfloor}x\bigr),
\qquad N\ge2,
\]
and established almost-everywhere convergence for every
$f\in L^p(X)$, $1<p<\infty$, whenever $h$ belongs to an
appropriate class of $c$-regularly varying functions.
This class includes the pure powers $h(t)=t^c, \ 1 < c < 4/3$, together with
perturbations such as
\[
h(t)=t^c(\log t)^a
\qquad\text{and}\qquad
h(t)=t^c\exp\bigl(a(\log t)^b\bigr),
\quad a\in\mathbb{R},\quad 0<b<1.
\]
	
	The integer-power case is covered by Nair's pointwise theorem \cite{Nair1993}.
	Here we establish variational estimates and pointwise convergence throughout
	the full range of noninteger exponents $c>1$. Our main result gives an
	$r$-variation estimate for every $r>2$. We use the homogeneous
	variation seminorm $\mathcal V^r$, defined in
	Section~\ref{sec:notation}.
	
	\begin{theorem}\label{t:maintheorem}
		Let $c>1$ be noninteger, let $1<p<\infty$, and let $r>2$.
		For every $\sigma$-finite measure-preserving system
		$(X,\mathcal B,\nu,T)$ and every $f\in L^p(X)$, we have
		\begin{equation}\label{eq:main-variation}
			\bigl\|\mathcal V^r(A_N^{c} f:N\in\mathbb N)\bigr\|_{L^p(X)}
			\lesssim_{c,p,r}\|f\|_{L^p(X)}.
		\end{equation}
		
	\end{theorem}
	
	 In	particular, standard considerations see for example \cite[Lemma 1]{FHK2007},  the estimate above gives the following pointwise convergence result,
	including the usual normalization of the prime averages.
	
	\begin{corollary}\label{cor:prime-convergence}
		Under the hypotheses of Theorem~\ref{t:maintheorem}, the averages
		$A_N^{c} f$ converge almost everywhere and in $L^p(X)$ as $N\to\infty$.
		The same conclusion holds for
		\[
		P_N^{c} f(x):=\frac1{\pi(N)}
		\sum_{p\le N,p \in \mathbb{P}}
		f\bigl(T^{\lfloor p^c\rfloor}x\bigr),\qquad N\ge2,
		\]
		and the two families have the same limit.
		Equivalently, if $p_j$ denotes the $j$-th prime, then
		\[
		\frac1M\sum_{j=1}^M f\bigl(T^{\lfloor p_j^c\rfloor}x\bigr)
		\]
		converges almost everywhere as $M\to\infty$.
	\end{corollary}

	We also consider the M\"obius-weighted averages
	\[
	M_N^{c} f(x):=\frac1N\sum_{n=1}^N
	\mu(n)f\bigl(T^{\lfloor n^c\rfloor}x\bigr).
	\]
	
    In an unpublished problem list originating from the conference "Pointwise Ergodic Theory \& Connections" in Bristol, Ter\"av\"ainen asked whether these
	averages converge almost everywhere to $0$ for every noninteger $c>1$ and $f \in L^2(X)$.
	The following result answers this question for $f\in L^p(X)$,
	$1<p<\infty$, and identifies the limit as well as a rate of convergence.

	\begin{theorem}[M\"obius-weighted convergence]\label{t:mobius-convergence}
		Let $c>1$ be noninteger and let $1<p<\infty$ and $A>1$. For every 
		$\sigma$-finite measure-preserving system $(X,\mathcal B,\nu,T)$
		and every $f\in L^p(X)$,
		\[
		\lim_{N\to\infty}\frac{\left(\log N\right)^A}{N}\sum_{n=1}^N
		\mu(n)f\bigl(T^{\lfloor n^c\rfloor}x\bigr)=0
		\]
		for almost every $x\in X$.
	\end{theorem}
	
	\subsection{Methods and organization}
	
	The main arithmetic input is the bound
	\begin{equation}\label{eq:intro-prime-comparison}
		\sup_{\xi\in\mathbb T}
		\left|\frac1N\sum_{n=1}^N
		(\Lambda(n)-1)e\bigl(\xi\lfloor n^c\rfloor\bigr)\right|
		\lesssim_{c,A}(\log N)^{-A}
		\qquad(A>0).
	\end{equation}

	On the integer shift, this estimate compares the multipliers of
	$A_N^{c}$ and $B_N^{c}$.
	The analysis in \cite{BDDS24} uses a circle method decomposition
	with a single major arc at zero and the bottleneck to extend the theorem therein to a larger set of parameters $c>1$ is the treatment of the minor arc which relies on the Vaughan identity and van der Corput's estimates. Here, we avoid the circle method entirely by using discorrelation estimates on short intervals to replace the von Mangoldt function with its corresponding Cram\'er approximant.
	
	  This is possible because the fractional-power phase can be approximated
	by a polynomial one of sufficiently large degree in a very effective manner in short intervals. Consequently, roughly speaking, one may replace the fractional phase with a piecewise polynomial one up to a negligible error.  This observation
	also appears in the work of Konieczny and M\"ullner
	\cite{CK24} and is crucial in the analysis therein as well.

    The short-interval uniformity estimates of Matom\"aki, Shao, Tao and Ter\"av\"ainen
	\cite{Matom_ki_2023} provide the (deep) arithmetic input. With the aid of fundamental lemma of sieve theory we replace the Cram\'er model by a divisor sum. The
	M\"obius version of the short-interval estimates provides the
	arithmetic input for Theorem~\ref{t:mobius-convergence}.
    
    We remark that short-interval uniformity estimates have also been employed in \cite{KoutsogiannisTsinas2026} for the study of multilinear ergodic averages along Hardy-field sequences evaluated at primes: after decomposing the prime-weighted averages into residue classes modulo $W$, they use short-interval Gowers uniformity estimates to compare averages weighted by the $W$-tricked von Mangoldt function with their unweighted counterparts in $L^2$, in a qualitative sense. 
	
	Section~\ref{sec:notation} records the notation, the Cram\'er model,
	the sieve estimate and the short-interval uniformity results.
	The technical lemmata concern floor removal, the analytic input coming from studying the averages $B_N^c$
	and approximation of the multiplier symbol corresponding to the averages  $B_N^c$ on certain arithmetic progressions. In addition to that,
	correlations of arithmetic weights with fractional phases on short intervals are essentially controlled by short-interval uniformity norms, see Lemma
	\ref{l:gowerscontrol}.  The subsequent sections treat
	the main variational theorem, and the final section treats the
	M\"obius-weighted averages.

	\section{Notation}\label{sec:notation}
	
	\subsection{General conventions}
	
	We write $\N=\{1,2,\ldots\}$,\;
	$\N_0=\N \cup \left\{ 0 \right\} $ and
	$[N]=\{1,\ldots,\lfloor N\rfloor\}$ for $N\ge1$.
	All sums over a real interval are understood to run over its integer
	points. For a finite nonempty set $E\subset\mathbb Z$, we put
	\[
	\mathbb E_{n\in E}f(n):=\frac1{|E|}\sum_{n\in E}f(n).
	\]
	For an interval $I$, the notation $|I|$ in a discrete average
	denotes $\#(I\cap\mathbb Z)$; when $I=(X,X+H]$, the length
	is $H$, and $|I|=H+O(1)$. The notation $J/d$ denotes the real
	interval $\{t:dt\in J\}$; a sum over $J/d$ is again over integers.
	
	We set
	\[
	e(t):=e^{2\pi i t},\qquad
	\{t\}:=t-\lfloor t\rfloor,\qquad
	\|t\|:=\min_{m\in\mathbb Z}|t-m|.
	\]
	The torus $\mathbb T=\mathbb R/\mathbb Z$ is equipped with
	normalized Haar measure. When a frequency in $\mathbb T$ occurs
	in an expression that is not periodic in that frequency, we use its
	representative in $[-1/2,1/2)$. 
	
 We write $F\lesssim_{A}G$, or
	$F=O_{A}(G)$, when $|F|\le C_{A}G$ for a
	constant depending only on the indicated parameters.
	The notation $F \sim_AG$ means that both corresponding
	inequalities hold. Dependence on fixed parameters such as $c$, $p$,
	the polynomial degree, and a fixed interval-length exponent may be
	suppressed. Constants are uniform in the averaging scale,
	frequencies, and weights unless stated otherwise.
	
	For a bounded interval $I\subset \R$, let $c_I$ denote its center.
	The Taylor polynomial of degree $k$ for $t\mapsto t^c$ at $c_I$ is going to be denoted by $P_{c_I}^k.$

	\subsection{Function spaces, Fourier transforms and variation seminorms}
	
	The space $\ell^p(\mathbb Z)$ is defined with counting measure.
	On a finite nonempty integer set $I$ we use normalized norms:
	\[
	\|w\|_{L^p(I)}:=
	\bigl(\mathbb E_{n\in I}|w(n)|^p\bigr)^{1/p}
	\quad(1\le p<\infty),\qquad
	\|w\|_{L^\infty(I)}:=\max_{n\in I}|w(n)|.
	\]
	For a measure-preserving system, $L^p(X)$ is taken with respect to
	$\nu$. 
	On $\mathbb Z$ we use the shift $Tx=x-1$.
	
	For $f\in\ell^1(\mathbb Z)$, our convention for the Fourier transform is
	\[
	\widehat f(\xi)=\sum_{x\in\mathbb Z}f(x)e(-x\xi),
	\qquad
	g^\vee(x)=\int_{\mathbb T}g(\xi)e(x\xi)\,d\xi.
	\]
	Convolution is given by
	$(f*g)(x)=\sum_{y\in\mathbb Z}f(x-y)g(y)$, and $\delta_m$
	denotes the unit point mass at $m$. On $\mathbb R$, the Fourier
	transform $\mathcal F_{\mathbb R}$ and its inverse use the same
	sign convention. We use
	$\mathrm M$ for the Hardy--Littlewood maximal operator on the
	ambient space.
	
	For an ordered set $E\subseteq\mathbb R$ and $1\le r<\infty$,
	define
	\[
	\mathcal V^r(a_t:t\in E)
	:=\sup_{\substack{L\ge1\\t_0<\cdots<t_L,\ t_j\in E}}
	\left(\sum_{j=1}^L|a_{t_j}-a_{t_{j-1}}|^r\right)^{1/r}.
	\]
	This seminorm is applied pointwise to families of functions. It is
	homogeneous: no term controlling the value at a single time is
	included in its definition. Since all averages are convex operators, there is no loss of generality in using the homogeneous formulation.
	
	\subsection{Arithmetic functions and Cram\'er approximants}
	
	We write $\mathbb P$ for the set of primes and
	\[ \pi(N)=\#\{p\in\mathbb P:p\le N\}.\]
    The von Mangoldt and
	M\"obius functions are
	\[
	\Lambda(n)=
	\begin{cases}
		\log p,&n=p^j,\ p\in\mathbb P,\ j\ge1,\\
		0,&\text{otherwise},
	\end{cases}
	\qquad
	\mu(n)=
	\begin{cases}
		(-1)^j,&n\text{ is a product of }j\text{ distinct primes},\\
		0,&n\text{ is not square-free},
	\end{cases}
	\]
	with $\mu(1)=1$. Euler's totient function is denoted by $\varphi$.

	For $w\ge2$, set
	\[
	W(w):=\prod_{\substack{p\le w\\p\in\mathbb P}}p,
	\qquad
	\Lambda_{\Cramer,w}(n)
	:=\frac{W(w)}{\varphi(W(w))}
	\mathbf1_{\gcd(n,W(w))=1}.
	\]
	
	\begin{definition}[Cram\'er approximant adapted to an interval]
		\label{def:interval-cramer}
		Given $J=(X,X+H]$ with $X\ge3$, put
		\[
		Q_J:=\exp\bigl((\log X)^{1/10}\bigr),\qquad
		W_J:=\prod_{\substack{p\le Q_J\\p\in\mathbb P}}p.
		\]
		The Cram\'er approximant adapted to $J$ is
		\[
		\Lambda_{\Cramer,Q_J}(n)
		:=\frac{W_J}{\varphi(W_J)}
		\mathbf1_{\gcd(n,W_J)=1}.
		\]
		In particular, the cutoff $Q_J$ is determined by the left endpoint
		$X$ of $J$, rather than by its length $H$.
	\end{definition}
	
	We use the following interval form of the fundamental lemma of
	sieve theory; see \cite[Lemma 4.2]{KMTT25}.
	
	\begin{lemma}[Fundamental lemma of sieve theory]\label{fund-lemma}
		Let $I$ be a finite interval of integers, and suppose that
		$2\le w\le y\le |I|^{1/10}$. Put $W=W(w)$. There exist real
		weights $\lambda_d^\pm\in[-1,1]$, supported on $1\le d\le y$,
		such that, for every integer $n$,
		\[
		\sum_{d\mid n}\lambda_d^-
		\le\frac{\varphi(W)}W\Lambda_{\Cramer,w}(n)
		\le\sum_{d\mid n}\lambda_d^+.
		\]
		Moreover,
		\begin{equation}\label{eq:expecationdivsum}
			\mathbb E_{n\in I}\sum_{d\mid n}\lambda_d^\pm
			=\frac{\varphi(W)}W
			\left(1+O\left(\exp\left(-\frac{\log y}{\log w}\right)\right)\right).
		\end{equation}
		In particular,
		\[
		\mathbb E_{n\in I}
		\left|\Lambda_{\Cramer,w}(n)
		-\frac W{\varphi(W)}\sum_{d\mid n}\lambda_d^\pm\right|
		\lesssim\exp\left(-\frac{\log y}{\log w}\right).
		\]
		The implied constants are absolute, and the estimates are uniform
		in the position of $I$.
	\end{lemma}
	
	\subsection{Short-interval polynomial uniformity}
	\label{subsec:short-uniformity}
	
	For $s\in\mathbb N$ and a finite nonempty interval $I$ of
	integers, define
	\[
	\|w\|_{u^s(I)}
	:=\sup_{\substack{P\in\mathbb R[t]\\\deg P\le s-1}}
	\left|\mathbb E_{n\in I}w(n)e(-P(n))\right|.
	\]
	We also write $u^s[N,N+H]$ when the interval is $[N,N+H]$.
	This polynomial-correlation seminorm measures correlation with
	polynomial phases; it is distinct from the full Gowers $U^s$ norm.
	The sign in the definition is immaterial.
	
	The polynomial-phase specialization of
	\cite[Theorem 1.1(i)--(ii)]{Matom_ki_2023} gives the following input.
	
	\begin{proposition}[Short-interval polynomial uniformity]
		\label{prop:short-uniformity}
		Fix $s\in\mathbb N$ and $0<\varepsilon<3/16$.
		Let $X\ge3$ and
		\[
		X^{5/8+\varepsilon}\le H\le X^{1-\varepsilon},
		\qquad J=(X,X+H].
		\]
		For every $A>0$,
		\[
		\|\Lambda-\Lambda_{\Cramer,Q_J}\|_{u^s(J)}
		+\|\mu\|_{u^s(J)}
		\lesssim_{A,s,\varepsilon}(\log X)^{-A}.
		\]
		The implied constant is independent of $X$ and $H$ and is uniform
		over all real coefficients of the polynomial phases.
	\end{proposition}
	
	In particular, for any fixed $\theta\in(5/8,1)$ and
	$H=\lfloor X^\theta\rfloor$, the same bound holds with a constant
	depending on $A,s,\theta$. Equivalently, for every real polynomial
	$P$ of degree at most $s-1$,
	\[
	\left|\sum_{n\in J}
	(\Lambda(n)-\Lambda_{\Cramer,Q_J}(n))e(P(n))\right|
	+\left|\sum_{n\in J}\mu(n)e(P(n))\right|
	\lesssim_{A,s,\theta}|J|(\log X)^{-A}.
	\]

	\section{Analytic input for fractional power orbits}
	We use the following refinement of van der Corput's estimate for polynomial-like phases; see \cite[Theorem 2.9]{GrahamKolesnik1991}.

	\begin{lemma}\label{l:vdc}
		Let \(k\geq 0\) be an integer, \(X\in\mathbb{N}\), and
		\(I\subset (X,2X]\) be an interval. Suppose that
		\(f\in C^{k+2}(I)\) and that there exists a constant
		\(G>0\) such that
		\[
		|f^{(r)}(x)|\sim GX^{-r},
		\qquad 1\leq r\leq k+2,
		\]
		uniformly for \(x\in I\). Then
		\[
		\left|\sum_{n\in I}e(f(n))\right|
		\lesssim
		G^{\frac{1}{2^{k+2}-2}}
		X^{1-\frac{k+2}{2^{k+2}-2}}
		+\frac{X}{G}.
		\]
	\end{lemma}
	
	The following estimate is used to dispose errors that have to do with the floor removal lemma procedure.
	\begin{lemma}\label{floorremovalerror}
		For every $c\in(1,\infty) \setminus \Z $ and  $N,M\in\mathbb{N}_{\ge 2}$ we have that there exists $\alpha \coloneqq \alpha(c)>0$ such that
		\begin{equation}\label{Lemerroroferror} 
			\sum_{n\le N} \min\left\{1,\frac{1}{M\|n^c\|}\right\}\lesssim  \log(N)  \left(N\frac{\log(M)}{M}+N^{1-\alpha} ( \log(M))^{0.5} +M\right)
		\end{equation}
	\end{lemma}
	\begin{proof}
		The proof is standard and the only discrepancy with the proof of  \cite[Lemma 3.11]{BDDS24} is the use of Lemma \ref{l:vdc} instead of the usual van der Corput lemma. In particular, it is clear that it suffices to control  \[  \sup_{1 \leq K \leq K' \leq 2K \leq N} \sum_{n \in [K,K']}	 \min\left\{1,\frac{1}{M\|n^c\|}\right\} \lesssim \sup_{1 \leq K \leq K' \leq 2K \leq N}	\sum_{m \in \Z} |a_m| \left| \sum_{K \leq n \leq K'}e(mn^c) \right| \] for some coefficients $|a_m| \lesssim \min \left\{\frac{\log M}{M}, \frac{M}{|m|^2}\right\}.$ By \cite[Lemma 2.1]{DDMSS25} we have that \[\begin{split}
			\sum_{m \in \Z} |a_m| \left| \sum_{K \leq n \leq K'}e(mn^c) \right|  &\lesssim K\frac{\log(M)}{M}+\sum_{0<|m| \leq K^{3c}} \min \left\{\frac{\log M}{M},\frac{M}{|m|^2}\right\} K^{1-\alpha}+ K\sum_{|m|>K^{3c}} \min \left\{\frac{\log M}{M},\frac{M}{|m|^2}\right\} \\ & \lesssim K\frac{\log(M)}{M}+ K^{1-\alpha} \left(\sum_{|m| \leq M \left(\log(M)\right)^{-0.5}} \frac{\log M}{M} +\sum_{ M \left(\log(M)\right)^{-0.5}<|m| \leq K^{3c} }\frac{M}{|m|^2} \right)+ \sum_{|m|>K^{3c}} \frac{MK}{|m|^2} \\ & \lesssim K\frac{\log(M)}{M}+K^{1-\alpha} ( \log(M))^{0.5}+M K^{1-3c}
		\end{split} \] Standard dyadic control of the initial sum yields \[ \sum_{n \leq N} \min\left\{1,\frac{1}{M\|n^c\|}\right\}\lesssim \log (N)   \sup_{1 \leq K \leq K' \leq 2K \leq N}	\sum_{K \leq n \leq K'} \min\left\{1,\frac{1}{M\|n^c\|}\right\} \lesssim \log(N)  \left(N\frac{\log(M)}{M}+N^{1-\alpha} ( \log(M))^{0.5} +M\right) \]
	\end{proof}

    At this point we are ready to state the main analytic input which we will utilize to evaluate the  correlation of a divisor sum with a fractional phase. The lemma below is probably well known to hold as the arguments to obtain it are rather standard  however we could not locate it in the literature.
	\begin{lemma}\label{l:analyticlemmacirclemethod} Let $c \in (1,\infty) \setminus \Z $ and $1 \leq d \leq \exp \left(\log(N)^{0.3}\right).$ Then there exists a $\delta_c>0$ such that  \[ \EE_{n \in [N/d]}e(\xi \floor{(dn)^c} )=\frac{d}{N} \int_{1}^{N/d} e(\xi (dt)^c) \d t+O(N^{-\delta_c})=\frac{1}{N} \int_{1}^Ne(\xi t^c) \d t+O(N^{-\delta_c})  \] holds uniformly in $\xi \in \TT.$ Namely the constant implied by the big-O notation does not depend on $\xi.$
		
	\end{lemma}
	\begin{proof}
		The second equality is rather straightforward to show and follows from the fact that $d$ satisfies $d\leq N^{o(1)}$. Fix a small parameter $\eps \in (0, \frac{1}{50}2^{-100\ceiling{c}^2})$.	If $|\xi| \leq N^{2\eps-c}$ then we have that \[\begin{split}
			&	\left| \sum_{n \in [N/d-1]}\left(\int_{n}^{n+1} \left(e(\xi \floor{(dn)^c})- e(\xi (dt)^c)\right) \d t \right) \right| \lesssim \sum_{n \in [N/d-1]} |\xi| \int_{n}^{n+1} | \floor{(dn)^c}-(dt)^c| \d t  \lesssim \\ &  \sum_{n \in [N/d-1]} |\xi| \max_{ t\in [n,n+1]}(dc(dt)^{c-1}+1 )\lesssim  \frac{N}{d}|\xi| +|\xi| \sum_{n \in [N/d-1]} d (dn)^{c-1+\eps} \lesssim \frac{N}{d}|\xi| + N^{\eps+c}|\xi|
		\end{split}\] 
        Hence, \[  \left|\EE_{n \in [N/d]}e(\xi \floor{(dn)^c} )-\frac{d}{N} \int_{1}^{N/d} e(\xi (dt)^c) \d t \right| \lesssim \frac{d}{N}+	|\xi|+dN^{\eps+c-1}|\xi|  \lesssim N^{-\delta_c}. \]

        On the other hand, if $|\xi| >N^{2\eps-c}$ by applying the van der Corput's lemma  since the modulus of the $\ceiling{c}$-derivative of the phase $\xi t^c$ is controlled from below by $|\xi t^{c-\ceiling{c}} | \gtrsim  N^{2\eps-\ceiling{c}}  $ gives \[ \left|\frac{1}{N} \int_{1}^N e(\xi t^c) \d t \right| \lesssim  \frac{1}{N \left( N^{2\eps-\ceiling{c}}  \right)^{\frac{1}{\ceiling{c}}}} \lesssim N^{-\frac{2\eps}{\ceiling{c}}}.   \]

        In addition to that, we employ the usual floor removing technology Fourier expanding $e(-\xi \left\{y\right\})$ as detailed in \cite[Page 29]{LPE}. To be more precise, for every $M\ge 2$, $\xi \in[-1/2,1/2)\setminus\{0\}$ and $y\in\mathbb{R}$ we have that
		\[
		e(-\xi \{y\})=\sum_{|m|\le M}\frac{1-e(-\xi )}{2\pi i (m+\xi )}e(my)+O\bigg(\min\bigg(1,\frac{1}{M\|y\|}\bigg)\bigg)\text{,}
		\]
		which combined with \eqref{Lemerroroferror} yields
		\[\begin{split}
			\EE_{n \in [N/d]}e(\xi \floor{(dn)^c})&=\sum_{|m| \leq N^{\frac{\eps}{100}}} \frac{1-e(-\xi)}{2 \pi i(m+\xi)} \EE_{n \in [N/d]} e\left((m+\xi)(dn)^c\right)+O(N^{-\eps_c}).
			\end{split} \]

        Furthermore, controlling the above exponential sum by the corresponding dyadic components conceding a $\log(N)$ term, namely: \[
		\left| \EE_{n \in [N/d]}  e\left((m+\xi)(dn)^c\right)\right|
		\lesssim
		\frac{d}{N}
		+\frac{d\log N}{N}\,
		\sup_{1\le M\le M'\le 2M\le N/d}
		\left|\sum_{M\le n\le M'} e\left((m+\xi)(dn)^c\right) \right|.
		\]   
        Lemma  \ref{l:vdc} applied for   $k=\ceiling{c}$  and triangle inequality allows us to bound each exponential sum, for every $|m|\leq N^{\eps/100}$, by \[\begin{split}
			\max_{|m|\leq N^{\eps/100}}\left| \sum_{M \leq n \leq M'}  e\left (\left(m+\xi\right)(dn)^c\right )\right| &\lesssim\max_{|m|\leq N^{\eps/100}}\min\left\{ \exp\left(\log(N)^{0.4}\right)\left(\left(|m+\xi|M^c\right)^{\frac{1}{2^{\ceiling{c}+2} -2}}M^{1-\frac{\ceiling{c}+2}{2^{\ceiling{c}+2}-2}}+ \frac{M}{|\xi| M^c} \right),M\right\} \\ &\lesssim N^{1-\eps_c}.  
		\end{split} \] 
Consequently, decreasing $\eps_c>0$ if necessary, \[\max_{|m|\leq N^{\eps/100}}\left|\EE_{n\in[N/d]}e\left((m+\xi)(dn)^c\right)\right|\lesssim N^{-\eps_c},\] and the proof is complete. \end{proof}
	
	The following unweighted variation estimate is contained in \cite[Corollary~1.6]{OKeeffe2024}. For completeness, we give a different standalone proof.
	\begin{proposition} \label{p:ncvar} For all $c \in (1,\infty) \setminus \Z$ and all $r>2$ and all $ p \in (1,\infty)$ we have that 
		\[  \left \| \mathcal{V}^r( \EE_{n \in [N]}f(x-\floor{n^c}) : N \in \N ) \right \|_{\ell^p(\Z)}  \lesssim \left \| f \right \|_{\ell^p(\Z)}. \] 
	\end{proposition}
	\begin{proof}
		Initially, standard control of the variational norm from the  long and the square function of the short variations \[ \mathcal{V}^r \left(B_N^{c} f: N \in \N\right) \lesssim \mathcal{V}^r\left(B_{N}^{c} f:N \in 2^{\N}  \right)+\left(\sum_{k \geq0} \mathcal{V}^2\left (B_N^{c} f:N \in [2^{k},2^{(k+1)}]\cap\N\right )^2\right)^{\frac{1}{2}}, \quad B_N^{c} f \coloneqq \EE_{n \in [N]}f(x-\floor{n^c}).\] 
        
        For the long variation, interpolating the strong $(2,2)$ bound provided by \ref{l:analyticlemmacirclemethod} against the strong $(1,1)$ and $(\infty,\infty)$ bounds allows us to control, for some $\kappa=\kappa(c,p)>0$, \[ \left \| \mathcal{V}^r \left( B_{N}^{c} f-\left(\widehat{f}(\xi) \eta(\xi) V_{N}^c(\xi)\right)^{\vee}: N \in 2^{\N}  \right) \right \|_{\ell^p(\Z)} \lesssim \sum_{N \in 2^{\N_0}} N^{-\kappa} \left \| f \right \|_{\ell^p(\Z)} \lesssim \left \| f \right \|_{\ell^p(\Z)}, \; \; V_N^c(\xi)= \frac{1}{N} \int_{0}^Ne(-\xi t^c ) \d t    \] where $\eta  \in C^{\infty}(\R)$ with $\ind{\left[-\frac{1}{10},\frac{1}{10}\right]} \leq \eta \leq \ind{\left[-\frac{1}{4},\frac{1}{4}\right]}$ and can be inserted in the estimate of Lemma \ref{l:analyticlemmacirclemethod} due to the estimate we obtained for the exponential sum when $|\xi| \geq N^{2\eps-c}.$ Usual Van der Corput estimates reveal the oscillatory and stationary behavior of $V_N^c$ namely that \begin{equation} \label{eq:vnestimates}
			V_N^c(\xi)= \begin{cases}
				O\left( \max\left\{1,N^c |\xi| \right\} ^{-\frac{1}{\ceiling{c}}}\right) \\ 1+O(|\xi|N^c)
			\end{cases}.
		\end{equation}

        The majorization
        \[
        \sup_{N>0}
        \left|\mathcal{F}_{\R}^{-1}
        \left[(\mathcal{F}_{\R}f)(\xi)\,\eta(\xi)\,V_N^c(\xi)\right](x)\right|
        \lesssim \mathrm{M}f(x),
        \]
       combined with the behavior of $V_N^c$, outlined above,  and Krause's remark \cite[Remark 3.134]{KrauseBook} and the
        Magyar--Stein--Wainger transference principle \cite[Proposition 2.1]{MSW02}, yields
       \[
        \left\|\mathcal{V}^r\left(B_N^{c} f:N\in 2^{\N} \right)\right\|_{\ell^p(\Z)}
        \lesssim \|f\|_{\ell^p(\Z)}.
        \]
   
        We now address the short variation component. The treatment is through the general principle developed in \cite[Theorem 2.39]{MSZ20}. For the application of this theorem adapted to our setting we use the Littlewood-Paley projections defined by means of their Fourier transform \[ \widehat{S_0f}(\xi)=\left(1-\psi(\xi)\right)\widehat{f}(\xi), \quad  \widehat{S_jf}(\xi)= \left(\psi(2^{(j-1)c}\xi)-\psi(2^{jc}\xi)\right)\widehat{f}(\xi), \quad j \geq 1 \]  where $\psi \in C^{\infty}(\R)$ is a bump with $ \ind{[-1/4,1/4]} \leq \psi \leq \ind{(-1/2,1/2)}.$ Set $S_j=0$ for $j<0$. Then it is clear that \[ \sum_{j  \geq 0} S_jf=f \] by working in frequency and summing the telescoping series. Furthermore, it is evident that \[ \left \| \left( \sum_{j \geq 0}|S_jf|^2 \right)^{\frac{1}{2}} \right \|_{\ell^p(\Z)} \lesssim \left \| f \right \|_{\ell^p(\Z)} \] by means of Khintchine's inequality and observing that the randomized multiplier is H\"{o}rmander-Mikhlin. 
        
        Before we verify the square function estimate, we first quickly show that \[ \left \| B_{N+h}^{c}-B_N^{c} \right \|_{\ell^p(\Z) \to \ell^p(\Z)} \lesssim \frac{h}{N} \]  for all $p\in[1,2]$ and $N,h\in\N$. We verify this for $p=1$ and $p=2$ and the interior follows from Riesz-Thorin interpolation. We write \[ \begin{split}
			B_{N+h}^{c} f-B_N^{c} f&= \frac{1}{N+h} \sum_{n=1}^{N+h}f(x-\floor{n^c})-\frac{1}{N} \sum_{n=1}^Nf(x-\floor{n^c})= \frac{-h}{N(N+h)} \sum_{n=1}^{N+h}f(x-\floor{n^c})+\frac{1}{N} \sum_{n=N+1}^{N+h}f(x-\floor{n^c}) \\ &  = \int_{\TT} \widehat{f}(\xi)e(\xi x) \left( \frac{-h}{N(N+h)} \sum_{n=1}^{N+h}e(-\xi \floor{n^c})+\frac{1}{N} \sum_{n=N+1}^{N+h} e(-\xi \floor{n^c})  \right) \d \xi 
		\end{split} \] hence $\left \| B_{N+h}^{c}-B_N^{c} \right \|_{\ell^1(\Z) \to \ell^1(\Z)} \lesssim \frac{h}{N}$ simply by the triangle inequality and $\left \| B_{N+h}^{c}-B_N^{c}\right \|_{\ell^2(\Z) \to \ell^2(\Z)} \lesssim \frac{h}{N}$ by bounding the corresponding Fourier multiplier symbol. At this point, we are ready to verify the conditions of \cite[Theorem 2.39]{MSZ20} and conclude that the square function of the short variations is $\ell^p(\Z)$ bounded as the corresponding maximal function is bounded on all $\ell^p$ for $p \in (1,\infty)$, where it is crucial that we can write $B_N^{c}$ as a convolution operator namely \[ B_N^{c} f(x)= \left(f* \frac{1}{N} \sum_{n=1}^N \delta_{\floor{n^c}} \right)(x). \] 
        
        Our first task is to verify the square function estimate  by noting that there is a $\sigma_c>0$ such that, for $j\in\Z$, $k,l\in\mathbb N_0$ with $k\geq l$, and $m\in\{0,\ldots,2^l-1\}$, \[ \begin{split}
			\left \| \left ( B_{2^k+2^{k-l}(m+1)}^{c}-B_{2^{k}+2^{k-l}m}^{c}\right)S_{j+k}f \right \|_2 & \lesssim \min \left\{ \left \| B_{2^k+2^{k-l}(m+1)}^{c}-B_{2^{k}+2^{k-l}m}^{c} \right \|_{\ell^2 \to \ell^2},2^{-|j|\sigma_c}  \right\} \left \| S_{j+k} f \right \|_{\ell^2(\Z)} \\ & \lesssim \min \left\{2^{-l},2^{-|j|\sigma_c}\right\} \left \| S_{j+k}f \right \|_{\ell^2(\Z)} 
		\end{split}  \] using Lemma \ref{l:analyticlemmacirclemethod}, \eqref{eq:vnestimates} and noting that if $S_{j+k}f \neq 0$ then $k \geq -j.$ Consequently, \[ \begin{split}
			\left \| \left( \sum_{k \geq l} \sum_{m=0}^{2^l-1}\left| \left(B_{2^k+2^{k-l}(m+1)}^{c}-B_{2^k+2^{k-l}m}^{c}\right) S_{j+k}f\right|^2 \right)^{\frac{1}{2}} \right \|_{\ell^2(\Z)} & \lesssim \left( \sum_{k \geq l} \sum_{m=0}^{2^l-1} \min \left\{2^{-2l},2^{-2|j| \sigma_c}\right\} \left \| S_{j+k} f \right \|_{\ell^2(\Z)}^2 \right)^{\frac{1}{2}} \\ &  \lesssim  \min \left\{2^{-\frac{l}{2}},2^{\frac{l}{2}-|j| \sigma_c}\right\} \left  \| f \right \|_{\ell^2(\Z)} \\ &  \lesssim  2^{-\frac{l}{2}} \min \left\{1, 2^{l-|j| \sigma_c}\right\} \left \| f \right \|_{\ell^2(\Z)}
		\end{split}\] Setting $a_{j,l}= \mathbf{C}_c \min \left\{1,2^{l-|j| \sigma_c}\right\}$ we can immediately calculate, for every $\rho>0$ and $\eps\in(0,\rho)$ \[ \sum_{l \geq 0} \sum_{j \in \Z} 2^{-\eps l} a_{j,l}^{\rho}= \mathbf{C}_c^{\rho}\sum_{l \geq 0} \sum_{j \in \Z} 2^{-\eps l} \min \left\{ 1, 2^{l\rho-|j| \sigma_c \rho}  \right\} \lesssim \sum_{l \geq 0} 2^{-\eps l} (l+1)  \lesssim  1.  \]
	\end{proof}
	
	\section{Controlling the correlations of an arithmetic weight with a fractional phase via short-interval $u^s$ norms} 
	The idea, also underlying in the arguments in \cite{CK24} is that if $I \subseteq [N,2N]$ is a short interval i.e  $|I| \lesssim N^{\theta}$ for some $\theta \in (0,1)$ then if $P_{c_I}^k$ is the $k$-th degree polynomial MVT reveals that  \[ \sup_{\xi \in \TT}\sup_{n \in I} \left|e(\xi n^c)-e(\xi P_{c_I}^k(n))  \right| \lesssim \sup_{n \in I} \left|n^c-P_{c_{I}}^k(n) \right| \lesssim  N^{c-1+\theta-(1-\theta)k}  \] which implies that on a short interval we may exchange the fractional phase $e(\xi n^c)$ with a polynomial of large enough degree phase up to a power-saving error. We make this heuristic precise in the following lemma, first fixing a small positive parameter $\eps'=\frac{1}{100}2^{-100\ceiling{c}^2}.$
	\begin{lemma} \label{l:gowerscontrol} 	Let $c >1$, an arithmetic weight $w:\Z\to\mathbb C$, and $I \subseteq [N,2N]$  with the additional property   $|I| \lesssim N^{\theta}$ for some $\theta \in (0,1).$ Then, there is $k \in \N$  and $\eta_{\theta,c}$ depending only on $\theta$ and $c$  such that \[ \left|\EE_{ n \in I}w(n)e(\xi \floor{n^c})\right| \lesssim \log(N) \left \| w \right \|_{u^{k+1}(I)}+ \left \| w \right \|_{L^{\infty}(I)}  \EE_{n \in I} \min \left\{1,\frac{1}{N^{\eps'}\left \| n^c\right \|}\right\}+\left \|w \right \|_{L^{\infty}(I)}\log(N)N^{-\eta_{\theta,c}}
		\] 
		
	\end{lemma}
	\begin{proof}
		Arguing as in the proof of Lemma \ref{l:analyticlemmacirclemethod} we write \[ \EE_{n \in I}w(n)e(\xi \floor{n^c})=\sum_{|m| \leq N^{\eps'}}\frac{1-e(-\xi)}{2 \pi i(m+\xi)}\EE_{n \in I }w(n)e\left((m+\xi)n^c\right)+O\left( \left \| w \right \|_{L^{\infty}(I)} \EE_{n \in I } \min \left\{1,\frac{1}{N^{\eps'} \left \| n^c \right \| }\right\} \right)\] and therefore by the calculation above and choosing $k$ large enough we have that there is $\eta_{\theta,c}>0$ such that we may conclude that \[ \EE_{n \in I}w(n)e\left(\left (m+\xi\right )n^c\right)=\EE_{n \in I}w(n)e\left (\left (m+\xi\right )P_{c_I}^k(n)\right )+O(\left \| w \right \|_{L^{\infty}(I)}N^{-\eta_{\theta,c}}). \] In turn, this implies that \[ \EE_{n \in I}w(n)e(\xi \floor{n^c})= \sum_{|m| \leq N^{\eps'}}\frac{1-e(-\xi)}{2 \pi i(m+\xi)} \EE_{n \in I}w(n)e\left (\left (m+\xi\right )P_{c_I}^k(n)\right )+O \left(\left \| w \right \|_{L^{\infty}(I)} \left[\EE_{n \in I } \min \left\{1,\frac{1}{N^{\eps'} \left \| n^c \right \| }\right\}  +\log(N)N^{-\eta_{\theta,c}}\right]\right). \] The desired inequality is immediate from the definition of the $u^{k+1}(I)$ norm and summing the harmonic series. 
	\end{proof}
	\begin{proposition} Let $c\in(1,\infty)\setminus\Z$. For every $A>0$, we have that
		\[
		\sup_{\xi\in\TT}\left|\EE_{n\in[N]}(\Lambda(n)-1)e\left(\xi\floor{n^c}\right)\right|
		\lesssim_{c,A}(\log N)^{-A}.
		\]
	\end{proposition}
	
	\begin{proof}Fix $\theta\in \left( \sqrt{\frac{5}{8}+\frac{1}{100}},1\right)$ and $\gamma\in(0,1)$, to be chosen sufficiently close to $1$. Initially, we have that for some $\mu_c>0$ \[ \begin{split}
			\left| \EE_{n \in [N]} \left(\Lambda(n)-1\right)e(\xi \floor{n^c}) \right| &\lesssim \frac{\log(N)}{N} \sup_{1 \leq K \leq K' \leq 2K \leq N} \left| \sum_{K \leq n \leq K'} \left(\Lambda(n)-1\right)e(\xi \floor{n^c}) \right|  \\ & \lesssim N^{-\mu_c}+\frac{\log(N)}{N} \sup_{\substack{1 \leq K \leq K' \leq 2K \leq N \\  |[K,K']| > K^{\theta}, K \geq N^{\gamma} }} \left| \sum_{K \leq n \leq K'} \left(\Lambda(n)-1\right)e(\xi \floor{n^c}) \right| 
			\end{split} \]  Partitioning each interval $[K,K']$ into subintervals, forming a collection which we are going to call $\mathcal{J}$ (depending on $K,K',\theta$ but we will omit this from the notation), each of which has length $ \sim (K'-K)^{\theta} \gtrsim K^{\theta^2} \geq K^{\frac{5}{8}+\frac{1}{100}}$ gives \[ \sum_{K \leq n \leq K'} \left( \Lambda(n)-1 \right)e(\xi \floor{n^c})=\sum_{J \in \mathcal{J}} \sum_{n \in J} \left(\Lambda(n)-1\right)e(\xi \floor{n^c}).  \] 
            
            Working directly on every short subinterval $J$ in the collection $\mathcal{J}$ we see that \[ \left| \sum_{J \in \mathcal{J}}\sum_{n \in J}\left(\Lambda(n)-1\right)e(\xi \floor{n^c}) \right| \leq \left|\sum_{J \in \mathcal{J}}\sum_{n \in J}\left(\Lambda(n)-\Lambda_{\Cramer,Q_J}(n)  \right)e(\xi \floor{n^c}) \right|+\left| \sum_{J \in \mathcal{J}}\sum_{n \in J} \left(\Lambda_{\Cramer,Q_J}(n)-1\right)e(\xi \floor{n^c}) \right| \] and due to Lemmata \ref{l:gowerscontrol},\ref{floorremovalerror} and Proposition \ref{prop:short-uniformity} the first term on the right hand side is controlled by \begin{equation} \label{eq:switchtocrammerest} \sum_{J \in \mathcal{J}}\left| \sum_{n \in J}\left(\Lambda(n)-\Lambda_{\Cramer,Q_J}(n)  \right)e(\xi \floor{n^c}) \right| \lesssim_A N \left(\log(N)\right)^{-A} 
		\end{equation}  
        
        For the other term we use the fundamental lemma of sieve theory \ref{fund-lemma} to switch from the Cram\'er approximant to a divisor sum. Specifically, we take $y=\exp\left(\left(\log(\inf J)\right)^{0.3}\right)$ and $w=Q_J$ as well as  \[ \left| \sum_{n \in J} \left(\Lambda_{\Cramer,Q_J}(n)- \frac{W_J}{\varphi(W_J)} \sum_{d \mid n} \lambda_{d}^{+} \right) \right| \lesssim |J| \exp \left(-\frac{\log\left(\inf J\right)^{0.3}}{\log(\inf J)^{0.1}}\right) \lesssim |J| \exp \left(-\mathbf{c}(\log(N))^{0.2}\right) \]
	for some small constant $\mathbf{c}>0.$ Having that said, we calculate \[ \frac{W_J}{\varphi(W_J)} \sum_{n \in J}\sum_{d \mid n}\lambda_d^+ e(\xi \floor{n^c})=\frac{W_J}{\varphi(W_J)} \sum_{d \leq y}\lambda_{d}^+ \sum_{n \in J/d} e(\xi \floor{(dn)^c}). \] 
    
    Using Lemma \ref{l:analyticlemmacirclemethod} at scales $\sup J$ and $\inf J$ and taking the difference shows that \[ \sum_{n \in J/d} e( \xi \floor{(dn)^c})=\frac{1}{d}\int_{J}e(\xi t^c) \d t +O(N^{1-\gamma  \delta_c}).  \] 
    In addition to this, by \eqref{eq:expecationdivsum} we have that \[ \EE_{n \in J} \sum_{d \mid n} \lambda_d^+=\sum_{d \leq y} \frac{\lambda_d^+}{d}+O\left(\frac{y}{|J|}\right)=\frac{\varphi(W_J)}{W_J}\left(1+O\left(\exp\left(-\mathbf{c}(\log(N))^{0.2}\right)\right)\right) \]  hence  \[\begin{split}
		\frac{W_J}{\varphi(W_J)} \sum_{n \in J}\sum_{d \mid n}\lambda_d^+ e(\xi \floor{n^c})& =\frac{W_J}{\varphi(W_J)} \sum_{d \leq y}\lambda_d^+\sum_{n \in J/d}e(\xi \floor{(dn)^c})=\frac{W_J}{\varphi(W_J)} \sum_{d \leq y  }\lambda_{d}^+ \left(\frac{1}{d} \int_{J}e(\xi t^c) \d t+O(N^{1-\delta_c \gamma})\right) \\ & =\frac{W_J}{\varphi(W_J)} \sum_{d \leq y  } \frac{\lambda_d^+}{d} \int_J e(\xi t^c) \d t +O\left( \frac{W_J}{\varphi(W_J)} y N^{1-\gamma\delta_c}\right) \\ &=  \int_J e( \xi t^c) \d t  +O \left(\frac{W_J}{\varphi(W_J)}y N^{1-\gamma \delta_c}+|J|\exp \left( -\mathbf{c}\left(\log(N)\right)^{0.2}\right) \right) .
	\end{split}\]  
where $\mathbf{c}$ is a small positive absolute constant.    
    Putting everything together, we learn \[ \begin{split}
		\left|\sum_{n \in J} \left( \Lambda_{\Cramer,Q_J}(n)-1 \right)e(\xi \floor{n^c})\right|  &\leq \left|\sum_{n \in J }\left( \Lambda_{\Cramer,Q_J}(n)-\frac{W_J}{\varphi(W_J)} \sum_{d \mid n} \lambda_d^+\right)e(\xi \floor{n^c})\right|+\left| \sum_{n \in J} \left( \frac{W_J}{\varphi(W_J)} \sum_{d \mid n} \lambda_d^+-1 \right)e(\xi \floor{n^c}) \right|  \\ & \lesssim |J| \exp(-\mathbf{c}(\log (N))^{0.2})+ (\log(N))^{0.1}\exp\left( \left(\log(N)\right)^{0.3}\right)N^{1-\gamma\delta_c} .
	\end{split}\]
	
	Summing in $\mathcal{J}$ and using \eqref{eq:switchtocrammerest} with $A+2$ in place of $A$ gives \[\begin{split}
		\left|\sum_{K \leq n \leq K'} \left( \Lambda(n)-1 \right)e(\xi \floor{n^c})\right| &\lesssim_A N(\log(N))^{-A-2}+|[K,K']|\exp(-\mathbf{c}(\log (N))^{0.2})\\&+K^{1-\theta^2} (\log(N))^{0.1}\exp\left( \left(\log(N)\right)^{0.3}\right)N^{1-\gamma\delta_c} \\ & \lesssim_A N(\log(N))^{-A-2}+N\exp(-\mathbf{c}(\log (N))^{0.2}) +(\log(N))^{0.1} \exp((\log (N))^{0.3})N^{1+1-\gamma \delta_c -\theta^2}
	\end{split} \]
	
	Since $\lim_{\theta,\gamma \to 1^{-}}(2-\gamma \delta_c-\theta^2)=1-\delta_c<1$ the appropriate choice of $\theta$ and $\gamma $ in the interval $ (0,1)$ is possible. Choose them so that $\gamma\delta_c+\theta^2>1$. Multiplying the preceding bound by $\log(N)/N$ and including the initial $N^{-\mu_c}$ error gives the required $O_A((\log(N))^{-A})$ bound. \end{proof}

		\section{Proof of the main theorem}
		\begin{proof}[Proof of Theorem \ref{t:maintheorem}] By  the Calder\'{o}n transference principle it suffices to work on the integer shift. As usual the variation is going to be controlled from the long and short variations. In particular, having fixed $p \in (1,\infty)$ and $r>2$ we choose another parameter $\tau \in (0,1)$ whose exact value will be decided later. Set $N_k=\floor{2^{k^{\tau}}}$ for $k\geq0$, so that $N_0=1$.
			\[ \mathcal{V}^r(A_N^{c} f: N \in \N) \lesssim \mathcal{V}^r\left(A_{N_k}^{c} f: k\geq0\right)+\left(\sum_{k \geq 0}\left( \mathcal{V}^2\left(A_N^{c} f: N \in [N_k,N_{k+1}]\cap\N\right)\right)^2  \right)^{\frac{1}{2}}. \] 
            
            The long variation is treated by comparison to the variation of the orbits along $\left\{\floor{n^c}\right\}_{n \in \N}.$ The short variation is treated in a similar fashion as in \cite[Page 11-12]{BDDS24}. In further detail, \[ \left \|   \mathcal{V}^r\left(A_{N_k}^{c} f: k\geq0 \right) \right \|_{\ell^p(\Z)} \lesssim \left(1+\sum_{k\geq1} \left(\log(N_k)\right)^{-\frac{10}{\tau}}\right) \left \| f \right \|_{\ell^p(\Z)}+\left \|   \mathcal{V}^r\left(B_{N_k}^{c} f: k\geq0 \right) \right \|_{\ell^p(\Z)} \lesssim \left \| f \right \|_{\ell^p(\Z)} \] where the passage to the final inequality was due to Proposition \ref{p:ncvar}. 
            
            Furthermore, by Young's inequality we have that \[ \begin{split}
				\left \| A_{M+1}^{c}-A_{M}^{c} \right \|_{\ell^p(\Z) \to \ell^p(\Z)} & \leq \left \| \frac{1}{M+1} \sum_{n \in [M+1]} \Lambda(n) \delta_{\floor{n^c}}-\frac{1}{M} \sum_{n \in [M]} \Lambda(n)\delta_{\floor{n^c}} \right \|_{\ell^1(\Z)} \\ &  \lesssim \frac{1}{M(M+1)} \sum_{n \in [M+1]} \Lambda(n)+ \frac{1}{M} \sum_{M<n \leq M+1} \Lambda(n)  \lesssim \frac{\log(M)}{M}
			\end{split} \]  Set $q=\min\{2,p\}$, and $\tau=\min\left\{1,\frac{1}{100q'}\right\}$. By Minkowski's inequality, \[\begin{split}
				\left\|\left(\sum_{k \geq 0}\left( \mathcal{V}^2\left(A_N^{c} f: N \in [N_k,N_{k+1}]\cap\N\right)\right)^2  \right)^{\frac{1}{2}}\right\|_{\ell^p(\Z)} & \leq \left(\sum_{k \geq 0} \left( \sum_{M=N_k}^{N_{k+1}-1} 	\left \| A_{M+1}^{c}-A_{M}^{c} \right \|_{\ell^p(\Z) \to \ell^p(\Z)}  \right)^q \right)^{\frac{1}{q}} \left \| f \right \|_{\ell^p(\Z)} \\ & \lesssim \left(1+\sum_{k \geq 1} k^{q(2\tau-1)} \right)^{\frac{1}{q}} \left \| f\right \|_{\ell^p(\Z)} \lesssim \| f \|_{\ell^p(\Z)}
			\end{split} \]	 and the proof is complete.
		\end{proof}
		\section{M\"obius-weighted convergence}
We now prove Theorem \ref{t:mobius-convergence}, which answers the following question of Joni Ter\"{a}v\"{a}inen. 
		\begin{problem}
Prove that for each noninteger $c > 1$, the following averages converge to $0$ almost everywhere
\[ M_N^{c} f(x)\coloneqq \mathbb{E}_{n\in[N]} \, \mu(n) f\bigl(T^{\lfloor n^c \rfloor}x\bigr)\]
for $f \in L^2(X)$, where $\mu$ is the M\"{o}bius function.
        \end{problem}
Before we give a proof, we remark that the original question contains a second part which asks for the widest possible range of noninteger $c>1$ with the property that the conclusion extends to $L^1(X).$

\begin{proof}[Proof of Theorem \ref{t:mobius-convergence}]
  By Calder\'{o}n's transference principle, it suffices to establish the inequalities \[  \begin{split} & \left\| \mathcal{V}^r\left( \left(\log N\right)^A M_N^cf: N \in \N \right) \right\|_{\ell^p(\Z)} \lesssim \left\| f \right\|_{\ell^p(\Z)} \\ & \left\| \left( \log N \right)^A M_N ^c f\right\|_{\ell^p(\Z)} \lesssim  \left( \log N \right)^{-A} \left\| f \right\|_{\ell^p(\Z)}  \end{split} \] for some $r>2,$ all $p \in (1,\infty)$ and $A>1$.  
  
  By noting that $ \| M_N^c \|_{\ell^{\infty}(\Z) \to \ell^{\infty}(\Z)} \leq 1$ and $\left\| M_N^c\right\|_{\ell^1(\Z) \to \ell^1(\Z)} \leq 1$ it suffices to work on $\ell^2(\Z)$ and obtain the corresponding $\left( \log N \right)^{-A}$ gain for all $A>1$ in order to conclude the second inequality.  
  Controlling  \[ \| M_N^{c} f \|_{\ell^2(\Z)} \leq   \left\| \EE_{n\in[N]}\mu(n)e(\xi \floor{n^c}) \right\|_{L^{\infty}(\TT)} \| f \|_{\ell^2(\Z)} \] allow us to concentrate on the corresponding exponential sum on the right hand side. 
  
  For $K\geq3$, partitioning $[K,K']$ into intervals of length $\sim K^{\theta}$, with fixed $5/8<\theta<1$, and bounding the final shorter interval trivially, Lemma \ref{l:gowerscontrol}, Proposition \ref{prop:short-uniformity} and Lemma \ref{floorremovalerror} give, for every $B>0$,
    \[
    \sup_{K\leq K'\leq2K}\left|\sum_{n\in[K,K']}\mu(n)e(\xi\floor{n^c})\right|
    \lesssim_B K(\log K)^{-B},
    \]
    uniformly in $\xi\in\TT$. 
    
    Taking $B=100A+2$ and treating $K\leq N^{1/2}$ trivially, we obtain
    \[\begin{split}
    \left|\sum_{n\in[N]}\mu(n)e(\xi\floor{n^c})\right|
    &\lesssim\log N\sup_{1\leq K\leq K'\leq2K\leq N}
    \left|\sum_{n\in[K,K']}\mu(n)e(\xi\floor{n^c})\right|\\
    &\lesssim_A\log N\left(N^{1/2}+N(\log N)^{-70A-2}\right)
    \lesssim_A N(\log N)^{-70A}
    \end{split}\]  which immediately implies that \[ \left|\EE_{n\in[N]}\mu(n) e( \xi \floor{n^c})\right| \lesssim (\log(N))^{-70A} \] and so $\| M_{N}^c \|_{\ell^2(\Z) \to \ell^2(\Z)} \lesssim \left( \log N \right)^{-70A}.$ 
    
    Controlling the variation as in the proof of Theorem \ref{t:maintheorem}  we choose a parameter $\tau \in (0,1)$ whose exact value will be decided later. Set $N_k=\floor{2^{k^{\tau}}}$ for $k\geq0$, so that $N_0=1$. 
			\[ \mathcal{V}^r( \left( \log N \right)^A M_N^{c} f: N \in \N) \lesssim \mathcal{V}^r\left(\left( \log N_k \right)^AM_{N_k}^{c} f: k\geq0\right)+\left(\sum_{k \geq 0}\left( \mathcal{V}^2\left(\left( \log N \right)^AM_N^{c} f: N \in [N_k,N_{k+1}]\cap\N\right)\right)^2  \right)^{\frac{1}{2}}. \] 
            
            The long variation is addressed in the same way as in the proof of Theorem \ref{t:maintheorem}, with the decay of $\left\| M_N^c \right\|_{\ell^p(\Z) \to \ell^p(\Z)}$ obtained above. 
            
            For the short variation we record the estimate \[ \left\| \left( \log (N+1) \right)^A M_{N+1}^c-(\log N)^A M_{N}^c \right\|_{\ell^p(\Z) \to \ell^p(\Z)} \lesssim \frac{\left( \log N \right)^A}{N} \] which follows from analogous considerations as the one for $A_N^c$ and hence  setting $q=\min\{2,p\}$, and $\tau=\min\left\{1,\frac{1}{100(A+1)q'}\right\}$ we get   \[\begin{split}
			&	\left\| \left(\left(\sum_{k \geq 0}\left( \mathcal{V}^2\left(\left(\log N \right)^A M_N^{c}f: N \in [N_k,N_{k+1}]\cap\N\right)\right)^2\right) \right)^{\frac{1}{2}} \right\|_{\ell^p(\Z)} \\ &  \leq \left(\sum_{k \geq 0} \left( \sum_{T=N_k}^{N_{k+1}-1} \left \| \left( \log(T+1)  \right)^A M_{T+1}^{c}- \left( \log(T) \right)^AM_{T}^{c} \right \|_{\ell^p(\Z) \to \ell^p(\Z)}  \right)^q \right)^{\frac{1}{q}} \left \| f \right \|_{\ell^p(\Z)} \lesssim \left(1+\sum_{k \geq 1}k^{q(\tau(A+1)-1)}   \right)^{\frac{1}{q}} \left \| f\right \|_{\ell^p(\Z)}  \\ & \lesssim \left\| f  \right\|_{\ell^p(\Z)}
			\end{split} \]	and the proof is concluded. 
\end{proof}
		\bibliography{fractionalprimesbib}
		\bibliographystyle{amsplain}
	\end{document}